\documentclass[10pt]{article}
\usepackage{amsmath,amssymb}
\usepackage{color}
\usepackage{mathrsfs}
\usepackage{mathptmx}
\usepackage{verbatim}
\usepackage{epsfig}
\usepackage{CJK}
\usepackage{bm}
\usepackage{amsfonts}
\usepackage{amsthm}
\input{epsf.sty}
\usepackage{epsfig}
\usepackage{setspace}
\def\<{\langle}
\def\>{\rangle}
\def\p{\partial}

\def\zd{\mathbb Z^d}

\def\e{\epsilon}
\def\be{\begin{equation}}
\def\ee{\end{equation}}
\def\ba{\begin{array}}
\def\ea{\end{array}}

\newtheorem{Theorem}{Theorem}

\newtheorem{Corollary}{Corollary}

\numberwithin{equation}{section}

\def\e{\epsilon}

\def\be{\begin{equation}}
\def\ee{\end{equation}}
\def\br{\begin{eqnarray}}
\def\er{\end{eqnarray}}
\def\p{\partial}

\title{KAM theorems and open problems for infinite dimensional Hamiltonian with short range
\thanks
{Supported by   NSFC11271076 and  NSFC11121101. }}
\author{$\mbox{Xiaoping YUAN}$\footnote{
E-mail:xpyuan@fudan.edu.cn }
\hspace{12pt}  \\
}
\begin{document}
 \maketitle

\begin{quote}
\small {\bf Abstract.} Introduce several KAM theorems for infinite dimensional Hamiltonian with short range and discuss the relationship between spectra of linearized operator and invariant tori. Especially, introduce a KAM theorem in the paper\cite{[Y]} published in CMP 2002, which shows that there are rich KAM tori for a class of Hamiltonian with short range and with linearized operator of pure point spectra. Here are also presented several open problems.
\end{quote}

\section{infinite dimensional KAM tori for infinite dimensional Hamiltonian with short range }

In 1986, Fr\"ohlich, Spencer and Wayne\cite{[F-S-W]} generalized Kolmogorov-Arnold-Moser (KAM) theorem to infinite dimensional Hamiltonian
\begin{equation}\label{1}H=\sum_{j\in\zd}\omega_jI_j+\e\sum_{\<i,j\>}f_{\,i,j\>} (I_i,I_j;\phi_i,\phi_j)\end{equation}
with symplectic structure $\sum_j dI_j\wedge d\phi_j$ where $f$ is of short range:
\be|f_{\,i,j\>} (I_i,I_j;\phi_i,\phi_j)|\le \text{const.}\;e^{-m|i-j|} \ee
for some $m>0$ or  it is of finite range:
\be f_{\<i,j\>}\equiv 0,\quad \text{unless}\;|i-j|\le \tilde \rho\ee
for some finite, positive range. Assume the perturbation $f_{\<i,j\>}$ is analytic on some domain and are $O(I^2)$. Also assume the frequencies
$\omega_i$ ($i\in\in\zd$) are weakly correlated, e.g., to be i.i.d. random variables with a smooth distribution such as
\be d\rho (\omega_i)=\bigg\{\begin{array}{ll}
\frac{2}{\sqrt{\pi}}\exp(-\omega_i^2),&\mbox{for any }\; \omega_i>0, \\
0,&\text{if} \;\omega_i\le 0.
\end{array}\ee
Considering initial conditions strongly localized in space: $I^0=(I^0_j\in\mathbb R_+:\;j\in\zd)$ with
\be  I_j^0=\exp(-|j|^{d+\alpha})\ee
with $\alpha>0$, they construct the following KAM theorem:
\begin{Theorem}\label{thm1}[Fr\"ohlich-Spencer-Wayne\cite{[F-S-W]},1986]

There exists $\e_0 > 0$  such that, for $\e<\e_0$, there is a set,
$\Omega(I^0)$ of ("nonresonant") frequencies, $\omega$, with Prob($\Omega(I^0)$) arbitrarily close
to one (depending on $\e$ and $\e_0$) such that if $\omega\in\Omega(I^0)$ then there is a
sequence, $(\tilde I_J)_{j\in\zd}$, of action variables with the properties that $|I_j^0-\tilde I_j^0|<\rho_j^0$
and that the trajectory of the canonical flow generated by the Hamiltonian,
H, of (\ref{1}) with initial conditions $(\tilde I,\tilde\phi)$ (for some $\tilde\phi_j\in[0, 2\pi), \forall\; j\in\zd)$  lies
on an infinite dimensional invariant torus, $\mathbb T(\tilde I)$.
\end{Theorem}

A similar result is in\cite{[V-B]} and\cite{[P1]}. This theorem is of wide range of application. For example, it can be used to the classical mechanics of a vibrating crystal lattice. A each site $j\in\zd$ of a crystal lattice, a  classical oscillator
(an atom or ion) is attached. Its configurations are described by vectors
$q_j\in\mathbb R^n$. The equations of motion of these oscillators are given by
\be\label{1.6} \ddot{q}_j=-({\omega^0_j})^2 q_j+\sum_if_{ij}(q_i-q_j)-\frac{\p U}{\p q_j}(q) \ee
where $f_{ij}$ and $\p U/\p q_j$ are of short range, e.g.,
\be f_{ij}=0=\frac{\p^2 U}{\p q_i\p q_j},\quad \text{unless}\; \;|i-j|=1.\ee
and the potential $U$ of the harmonic forces
is assumed to be analytic in the position variables
and $0\le U\le O(||q||_{\ell^2}^4)$. The matrices $(\omega^0_j)^2$  and $f_{ij}$ are independent, identically distributed (i.i.d.) $n\times n$
random matrices with smooth distributions of fast decay at infinity which
have their support on positive  matrices. The most interesting case for
crystal physics corresponds to $\omega^0_j\equiv 0$.
Following \cite{[F-S-W]}, the linear part of (\ref{1.6}) may be rewritten in more compact notation as follows:
\be\label{1.8}\ddot{q}_j=-(\Omega^2q)_j ,\; j\in\zd\ee
where $\Omega^2$ is the Jacobi matrix given by
\be (\Omega)_{ij}=\left\{\begin{array}{ll} -f_{ij}, &  \text{if}\;|i-j|=1\\ (\omega_i^0)^2+\sum_{k:\;|k-i|=1}f_{ki},&
\text{if}\; j=i\\0,&\text{otherwise.}\end{array}\right.\ee
If we define $-\triangle_f$  to be the off-diagonal part of $\Omega^2$
and
\be v(j)\equiv (\omega^0)^2+\sum_{k:\;|k-i|=1}f_{ki}\ee
then
\be \Omega^2=-\triangle_f+v.\ee
Equation (\ref{1.6}) may be rewritten as
\be\ddot q=(-\triangle_f+v)q-\frac{\p U}{\p q}(q).\ee
When $f\equiv 1$, it is the usual discrete nonlinear wave equation.  If $\omega^0_j$ has a smooth distribution of sufficiently large disorder, the
entire spectrum of $\Omega^2$ is pure point, consisting of exponentially decaying
eigenfunctions, $q^i$. In this case, (\ref{1.8}) has a complete, orthonormal system of periodic solutions
\be\label{1.13}\cos(\omega_i t)q_j^i,\quad \sin(\omega_i t)q_j^i,\quad i\in\zd\ee
with period $2\pi/\omega_i$. Every solution, $q(t)$, of (\ref{1.8}) is a linear combination of the solutions (\ref{1.13}):
\be q_j(t)=\sum_{i\in\zd}\left[Q_i\cos(\omega_i t)+\frac{P_i}{\omega_i}\sin(\omega_i t)  \right]q_j^i \ee
for some finite, real coefficients $Q_i,P_i$. The action $I_j$ can be defined as
\be I_j=Q_j^2+\left( \frac{P_j}{\omega_j}\right)^2.\ee
Define a torus $\mathbb T^\infty(\e)$ by
\be\mathbb T^\infty(\e)=\{I^0_j\in\mathbb R: I^0_j\equiv 2\e_j,\;j\in\zd.\} \ee
Theorem \ref{thm1} shows that for ``most" $\omega$ there is an invariant torus of (\ref{1.6}) around $\mathbb T^\infty(\e)$ if
\be\label{1.17} \e_j:=Q_j^2+\left( \frac{P_j}{\omega_j}\right)^2\thickapprox\exp(-|j|^{d+\alpha}).\ee This invariant torus carries quasi-periodic or almost periodic motion.
Theorem \ref{thm1} can also be applied to the discrete nonlinear Schr\"odinger equation with a random potential
has the form
\be\label{1.18} {\bf i}\frac{\p}{\p t}\psi_t(x)=[(-\triangle+v)\psi_t](x)+\lambda\psi_t(x)\sum_{y\in\zd}V(|x-y|)|\psi_t(y)|^2.\ee
See Section 4.2 in \cite{[F-S-W]} for details. Theorem \ref{thm1} is generalized to partial differential equations. See \cite{[Bo2]} and \cite{[P2]}.

The initial conditions (\ref{1.17}) implies that the action (or the amplitude) decays too fast such that the infinite dimensional KAM tori look like to be of finite dimension. Weaker restrictions were later suggested in P\"oschel\cite{[P1]}:
\[\e_j=I_j^0\approx \exp(-(\log|j|)^{1+\alpha})\]
with some $\alpha>0$. Kuksin\cite{[K2]} and Bourgain\cite{[Bo1]} proposed the following open problem:

\vspace{5pt}

\noindent{\bf Problem 1}. {\it Is there any infinite dimensional KAM tori for the initial conditions with polynomial decay
\[\e_j=I_j^0\approx |j|^{-M} \]
with some large $M>0$?}

\section{finite dimensional KAM tori for infinite dimensional Hamiltonian with short range }

The study of the finite dimensional KAM tori for infinite dimensional Hamiltonian defined by partial differential equations was initiated by Kuksin in the end of 1980's. There are now a lot of work on this field. Here we focus on the Hamiltonian of short range instead of PDEs.

In order to state the related results, we need introduce some notations. Denote by $\ell$ the Hilbert space
$$\ell=\{u=(u_n^1,u_n^2)_{n\in \Bbb Z}:(u_n^1,u^2_n)\in \Bbb C^2,||u||:=\sum_{n\in \Bbb Z}(|u_n^1|^2+|u_n^2|^2)
e^{|n|/a}<\infty \}$$
where the inner product $\<\cdot,\cdot\>_{\ell}$ in $\ell$  is defined as follows:
$$\<u,v\>_\ell:=\sum_{n\in \mathbb Z}(u_n\cdot v_n) e^{|n|/a},\quad \text{some}\; a>0.$$
Set
$$\mathcal D_0=\{(I,\theta,u)\in \Bbb C^N\times \mathbb C^N/2\pi\mathbb Z^n\times\ell:|I| <s_0^{2/3},|Im\theta|< \delta_0,||u||< s_0^{1/3}\},$$
where $s_0$ and $\delta_0$ are positive constants. Following \cite{[Y]}, consider the Hamiltonian of the form

\be\label{2.1} H(I,\theta,u;\xi)=\sum_{j=1}^N\omega_j(\xi) I_j+\sum_{n\in \Bbb Z}\frac12\beta|u_n|^2+\epsilon P(\theta,I,u;\xi)
+P^2+ P^3,\ee with symplectic structure
\[\sum_{j=1}^N dI_j\wedge d\theta_j+\sum_{j\in\mathbb Z}du_j^1\wedge du_j^2,\]
where $H$ satisfies the following conditions

(A1). Hamiltonian  $H$ is analytic in $\mathcal D_0 \times \Pi$  and real for real arguments, where $\Pi$ is a parameter set of positive Lebesgue measure.

(A2). On the domain $\mathcal D_0$ the Hamiltonian $H$ is  of short range:

\be \label{2.2}P^2(u;\xi)=\epsilon \sum_{n\in \Bbb Z}W(u_{n+1}^1-u_n^1)\ee
\be P^3(I,u;\xi)=\sum_{n\in \Bbb Z}O(|u_n|^3)+O(|I|^2)\ee
\be \frac{\partial P(I,\theta,u;\xi)}{\partial u_j}\equiv 0
\quad\text{if}\quad |j| > 1\ee
 where
$u=(u_n)_{n\in \Bbb Z}$,
$\partial P/\partial u_j =(\partial P/\partial u_j^1,\partial P/\partial u_j^2 )$ and $u_j =(u_j^1,u_j^2)$, and
\be W(x)=O(x^3).\ee

(A3).(Non-degenerate) There is a constant $\delta_a>0$ such that
$$|\frac{\partial \omega}{\partial \xi}|\ge \delta_a \quad \text{on}\quad \Pi.$$

(A4). $$|P|\le K_1,\quad ||\nabla_u P||\le K_1 s_0^{-1/3}$$
where $K_1$ is a positive constant,$\nabla$ is the gradient  with respect to  the usual inner product $\<\cdot,\cdot\>$
in the usual square-summable space $\ell^2$.

\begin{Theorem}\label{thm2}[Yuan\cite{[Y]},2002]  Suppose that the Hamiltonian (\ref{2.1})
satisfies conditions (A1--4).Then, for a given $\gamma>0$ there is a small constant $\epsilon^*=\epsilon^*(\Pi,N,\gamma,\delta_a)>0$
such that, if $0<\epsilon<\epsilon^*$,then there is a Cantor set $\Pi_\infty\subset \Pi$ with $\text{meas}\Pi_\infty\ge (\text{meas}\Pi)
\cdot(1-O(\gamma))$,an analytic family of torus embedding $\Psi_\infty:\mathbb T^N\times \Pi\to\mathcal P$,and a map $\omega_\infty:\Pi\to\mathbb R^N$,such that for each $\xi\in \Pi_\infty$,the map $\Psi_\infty$ restricted to $\mathbb T^N\times \{\xi\}$ is an analytic embedding of rotational torus with frequencies $\omega_\infty(\xi)$ for  the Hamiltonian $H$.
\end{Theorem}

The main aim to construct this theorem is to find so-called quasi-periodic breathers of the
networks of weakly coupled oscillators: \be\label{2.6}\frac{d^2 x_n^2}{dt^2}+V'(x_n)=\epsilon W'(x_{n+1}-x_n)-\epsilon W'(x_n-x_{n-1}),
\quad n\in \mathbb Z.\ee
where $V$ is the local potential with  $ V'(0)=0,V''(0)=\beta$,( $\beta>0$),
and $W$ is the coupling potential.This equation has been deeply investigated by
some authors. See \cite{[A3]}, \cite{[M-A]}, \cite{[F-G]}, \cite{[Ka1]} and \cite{[Ka2]} for example.
In the classical case breathers are time-periodic and spatially localized solutions
of the equations of motion. Aubry \cite{[A1]}-\cite{[A-A]}
posed the well-known concept of anti-integrability or anti-continuation by which the existence of  the breathers for
 some equations of motion
 can be proven. While searching periodic breathers, small divisor problem does not appear.
Thus, in 1994, the existence of periodic breathers in a wide class of models (including (\ref{2.6}) was proved  by
MacKay and Aubry \cite{[M-A]}, using the anti-integrability. Naturally, one hopes to investigate so-called quasi-periodic breathers, i.e., the solutions which are quasi-periodic in time and exponential decay (localization) in space.  When searching quasi-periodic breathers, here indeed appears the small divisor problem. Thus, Aubry\cite{[A2]} remarked that the existence of the quasi-periodic breathers  is  an open problem which should relate
the concepts of the KAM theory and of anti-integrability.  By using Theorem \ref{thm2}, the open problem was solved with the assumption that  the linearized system has no continuous spectrum.

\begin{Theorem}\label{thm3}[Yuan\cite{[Y]},2002]
If $V(0)=V'(0)=0,V''(0)>0$,$W=O(|x|^3)$, and there is $k\ge 3$ such that $V^{(3)}(0)=\cdots=V^{(k-1)}(0)=0$,
$V^{(k)}(0)\neq 0$, then (\ref{2.6}) has ``rich''  quasi-periodic breathers
of small amplitude when $\epsilon$ is sufficiently small.
\end{Theorem}

\vspace{5pt}

\noindent{\bf Remark 1.} Usually the breathers decay exponentially. Here the obtained quasi-periodic breathers decay super-exponentially:
\be |x_n|\le \text{const.}\; e^{-|n|/a}\exp\left(\frac13(\frac43)^{|n|} \right). \ee
This may be a new phenomenon, which was recently observed by Kastner\cite{[Ka1]}\cite{[Ka2]}, in numeric experiment.

\vspace{5pt}

\noindent{\bf Remark 2.} Since the potential $V$ in (\ref{2.6}) is independent of the site $n$, the normal frequency
\be \lambda_n=V^{\prime\prime}(0)=\beta,\ee
which implies that, in this case, the frequencies $\lambda_n$'s ($n\in\mathbb Z$) are the most ``dense", i.e., they are condensed at a point $\beta$. The main aim of Theorem \ref{thm2} is to prove the existence of quasi-periodic breathers for (\ref{2.6}). Therefore, we take $\lambda_n=\beta$ in Theorem \ref{thm2}. Actually, it is not necessary to impose any conditions on the frequencies $\lambda_n$, especially not necessary to assume that all frequencies $\lambda_n$ are the same: $\lambda_n=\beta$. In fact, the scheme of KAM iteration in the proof of Theorem \ref{thm2} is the same as that devised by Kuksin\cite{[K1]} for partial differential equations. In the both cases, one encounters the following small divisor problem:
\be \<k,\omega\>\neq 0, \quad \forall\;0\neq k\in\mathbb Z^N, \ee
\be  \<k,\omega\>+\lambda_n\neq 0, \quad \forall\;k\in\mathbb Z^N, n\in\zd,\ee
\be \<k,\omega\>+\lambda_n-\lambda_m\neq0, \quad \forall\;k\in\mathbb Z^N, m\neq n\in\zd. \ee
They imply the infinite number of small divisor conditions.  However, the number must be finite in each KAM iteration step. In the KAM theory for PDEs, the finiteness of the number of the small divisors can be fulfilled by the growth of the frequencies $\lambda_n$. For example, for nonlinear Schr\"odinger equation,
\be \lambda_n=n^2,\quad |\lambda_m-\lambda_n|\ge m+n,\quad \forall m,n\in\mathbb Z_+, m\neq n.\ee
First, by the analyticity of the system itself, there is $\mathcal K=\mathcal K_l>0$ at the $l$-step of KAM iteration such that
$|k|\le \mathcal K$. Moreover, when $m$ or $n$ is larger than $C\mathcal K$ with $C=1+\sup|\omega(\xi)|$,
\be|\<k,\omega\>+\lambda_n|\ge \mathcal K,\quad |\<k,\omega\>+\lambda_n-\lambda_m|\ge \mathcal K,\ee
which are not small.
Thus, the number of the small divisors is bounded by
\be \mathcal K^{N+2}.\ee
From here one sees that the interesting case for Hamiltonian with short range is to assume $\lambda_n$'s are bounded, especially, $\lambda_n=\beta$ for all $n$. In \cite{[Y]}, it was observed that {\it the finiteness of the number of the small divisors can be fulfilled via the short range}. See (3.25), (3.26) and (3.27) in \cite{[Y]}. In addition, the short range condition (\ref{2.2}) can be replaced by  more general one
\be \label{bu1} \left| \frac{\p^2 P^2}{\p u_n\p u_m}\right|\le\e \exp({-C|m-n|}),\ee
since  the part of $P^2$ which involves large $m$ and $|n|\le \mathcal K$ can be discarded at the $l$-step KAM iteration.
 Therefore, Theorem \ref{thm2} can be rewritten as follows, not necessary to add  more argument:

\begin{Theorem}\label{thm4}  Suppose that the Hamiltonian
\be H(I,\theta,u;\xi)=\sum_{j=1}^N\omega_j(\xi) I_j+\sum_{n\in \zd}\frac12\lambda_n|u_n|^2+\epsilon P(\theta,I,u;\xi)
+P^2+ P^3 \ee
satisfies conditions (A1--4), where (\ref{2.2}) is replaced by (\ref{bu1}) and $\mathbb Z$ replaced by $\zd$. Then, for a given $\gamma>0$ there is a small constant $\epsilon^*=\epsilon^*(\Pi,N,\gamma,\delta_a)>0$
such that, if $0<\epsilon<\epsilon^*$,then there is a Cantor set $\Pi_\infty\subset \Pi$ with $\text{meas}\Pi_\infty\ge (\text{meas}\Pi)
\cdot(1-O(\gamma))$,an analytic family of torus embedding $\Psi_\infty:\mathbb T^N\times \Pi\to\mathcal P$,and a map $\omega_\infty:\Pi\to\mathbb R^N$,such that for each $\xi\in \Pi_\infty$,the map $\Psi_\infty$ restricted to $\mathbb T^N\times \{\xi\}$ is an analytic embedding of rotational torus with frequencies $\omega_\infty(\xi)$ for  the Hamiltonian $H$.
\end{Theorem}

Here is no any conditions imposed on the normal frequencies $\lambda_n$'s, except that which are of pure point spectra and are non-zero. Therefore, these frequencies can densely distributed on an finite interval. For example, they can be all rational number in the interval $[1,2]$. Thus, by Theorem \ref{thm4}, we have

\begin{Corollary}\label{thm5} Assume $\Omega^2=-\triangle_f+v$ is of pure point spectra. Then there are rich quasi-periodic breathers for discrete nonlinear wave equation (\ref{1.6})
\be\ddot{q}_j=((-\triangle_f+v)q)_j-\frac{\p U}{\p q_j}(q) ,\; j\in\zd\ee
and discrete nonlinear Schr\"odinger equation (\ref{1.18}).
\be{\bf i}\dot q_n=-(q_{n+1}+q_{n-1})+v(n)q_n+\lambda q_n\sum_{m\in\zd}V(|n-m|)|q_m|^2,\; n\in\zd\ee
where $U, V$ are of short range.
\end{Corollary}

\vspace{5pt}

\noindent{\bf Remark 3.} We have seen that the pure point spectra property is of key importance for the existence of quasi-periodic breathers for Hamiltonian with short range. This property is fulfilled by assume $W(x)=O(|x|^3)$ in Theorem \ref{thm3}. If $W(x)=\frac12 x^2+\cdots$, the linearized system of (\ref{2.6}) is of continuous spectra which is the whole of the interval $[\beta,\beta+1]$. When searching periodic breathers, there is no small divisors. Thus, in Aubry-MacKay's work\cite{[M-A]}, the coupled potential $W$ may be $O(|x|^2)$, which involving continuous spectra. If searching quasi-periodic solution of some symmetry, small divisors can also be avoided. Consequently, in 2002, Bambusi and Vella \cite{[B-V]} found many (but with measure $0$) quasi-periodic breathers of some symmetry for $W(x)=c x^2+O(x^3)$.

\vspace{5pt}

\noindent{\bf Problem 2.} {\it Assume that the linearized operator ($-\triangle_f+v$, for example,) of a hamiltonian system is of pure point spectra, part of which is dense in some finite interval, and assume that the nonlinear part of the Hamiltonian is {\bf not} of short rang (in partial differential equations defined in the whole space, for example). Are there any KAM tori for the system?}

This problem is apparently open.
\vspace{5pt}

\noindent{\bf Problem 3.} {\it Assume that the linearized operator ($-\triangle_f+v$, for example,) of a hamiltonian system contains continuous spectra which consist of an interval. Are there any KAM tori for the system? even if the system may be of short range.}

\vspace{4pt}
In general, the answer trends to be no.

\vspace{5pt}

\noindent{\bf Remark 4.} In (\ref{2.6}), the Hamiltonian of the coupled potential can be written
\be W(\{x_n\})=\sum_{n\in\mathbb Z }O(|x_{n+1}-x_n|^3).\ee
In 2007, Geng and Yi\cite{[G-Y]} replaced it by
\be W(\{x_n\})=\sum_{m\in\mathbb Z,m\neq n}e^{-|n-m|^\alpha}|x_n-x_m|^2,\ee
and replaced $V(x)$ by $V_n(x)$ with
\be V_n^{\prime\prime}(0):=\beta^2_n, \quad |\beta_n-\beta_m|\ge \gamma.\ee
This implies that the linearized operator is of pure point spectra and the normal frequencies are unbounded:
\be\lambda_n=\beta^2_n\to\infty.\ee



 \end{document}